\documentstyle[12pt,amssymb,graphicx]{article}
\begin{document}
\baselineskip=13pt
\begin{center}
{\Large\bf Irreducibility of spatial graphs}
\end{center}
\vspace{5mm}
\centerline{\large\it Dedicated to Professor Shin'ichi Suzuki for his 60th birthday}
\vspace{5mm}
\centerline{\large Kouki Taniyama}

\vspace{5mm}{\bf Abstract}\ \  A graph embedded in the 3-sphere is called irreducible if it is
non-splittable and for any 2-sphere embedded in the 3-sphere that intersects the graph at one point
the graph is contained in one of the 3-balls bounded by the 2-sphere. We show that irreducibility is
preserved under certain deformations of embedded graphs. We show that certain embedded graphs
are irreducible.
\vfill

\footnoterule



\newpage
\noindent{\large{\bf\S1. Introduction}}

Throughout this paper we work in the piecewise linear category. We consider a graph
as a topological space in the usual way. Let
$G$ be a finite graph embedded in the 3-sphere $S^3$. A
compact surface $F$ embedded in $S^3$ is called {\it good for $G$} if $\partial F$ is contained in
$G$, ${\rm int}F\cap G$ contains at most finitely many points, and for each point $x\in {\rm int}F\cap
G$ a neighbourhood of
$x$ looks like Fig. 1.1 for some positive integers $p$ and $q$. Note that $x$ is not necessarily a
vertex of $G$. Namely we allow the case that $p=q=1$ and $x$ is an interior point of an edge of $G$.

\begin{center}
\includegraphics[width=2.5in]{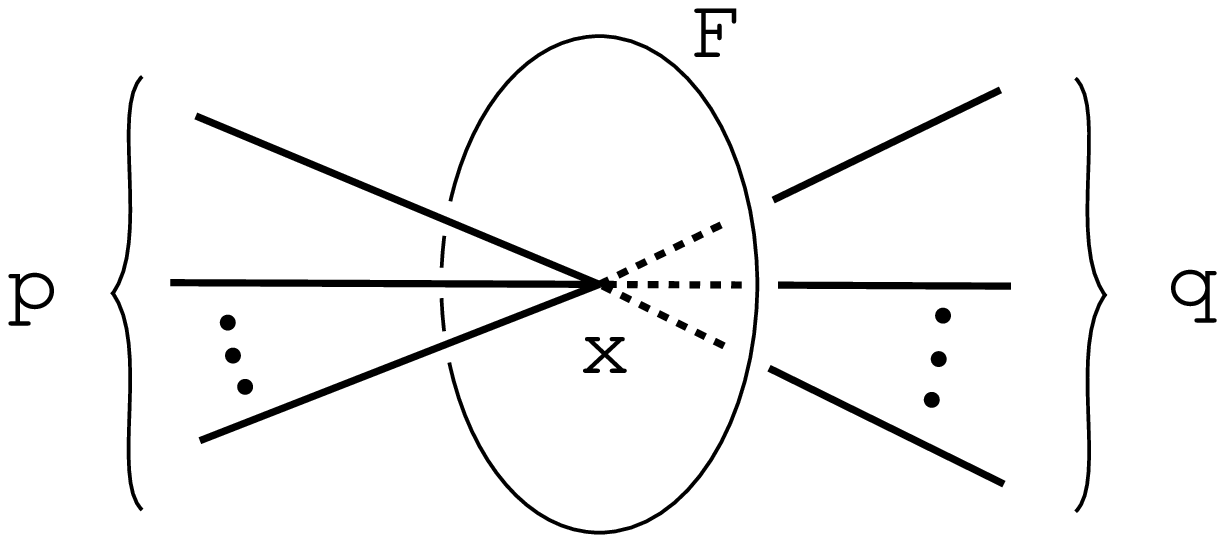}

Fig. 1.1
\end{center}

A 2-sphere $S$ embedded in $S^3$ is called a {\it separating sphere of $G$} if it is disjoint from
$G$ and each component of $S^3-S$ has intersection with $G$. A 2-sphere $S$ embedded in $S^3$ is
called a {\it cutting sphere of $G$} if it is good for $G$ and intersects $G$ at one point. Note
that if $S$ is a cutting sphere of $G$ then each component of
$S^3-S$ has intersection with $G$. A graph $G$ embedded in $S^3$ is called {\it non-splittable} if
there are no separating spheres of $G$. A graph $G$ embedded in $S^3$ is called {\it irreducible} if
there are no separating spheres and cutting spheres of $G$. Thus $G$ is non-splittable if $G$
is irreducible. The purpose of this paper is to show that certain graphs embedded in $S^3$ are
irreducible. Let
$D$ be a 2-disk embedded in
$S^3$ that is good for
$G$. Let
$G'$ be a graph embedded in $S^3$ obtained from $G$ by contracting $D$ to a point. Namely we
consider the quotient space $S^3/D$ and let $G'$ be the image of $G$ under the quotient map
$S^3\rightarrow S^3/D$. It is easy to see that $S^3/D$ is again
homeomorphic to $S^3$. By identifying $S^3/D$ with $S^3$ we may suppose that $G'$ is embedded in
$S^3$. Therefore $G'\subset S^3$
is well-defined up to self-homeomorphism of $S^3$.

\vspace{5mm}\noindent{\bf Theorem.} {\it Let $G$ be a graph embedded in $S^3$ and $D$ a 2-disk
embedded in $S^3$ that is good for $G$. Suppose that ${\rm int}D\cap G$ is not an empty set or
$\partial D\cap{\rm cl}(G-\partial D)$ is not a singleton where ${\rm cl}$ denotes the closure. Let
$G'$ be a graph embedded in
$S^3$ obtained from
$G$ by contracting $D$ to a point. If $G'$ is
irreducible then
$G$ is irreducible.}

\vspace{5mm}
In the following we show some applications of Theorem. 

The link illustrated in Fig. 1.2 (a) is a
Borromean ring. By contracting the obvious disk we have a graph illustrated in Fig. 1.2 (b). Again by
contracting another obvious disk we have a graph illustrated in Fig. 1.2 (c). Since this graph has
no loops and cut-vertices in the sense of abstract graph we have that this graph is irreducible.
Then by Theorem we have that the graph illustrated in Fig. 1.2 (b) is also irreducible. Then
again by Theorem we have that the Borromean ring is also irreducible, hence it is non-splittable.
This is an elementary geometric proof of the nontriviality of a Borromean ring.

\begin{center}
\includegraphics[width=5in]{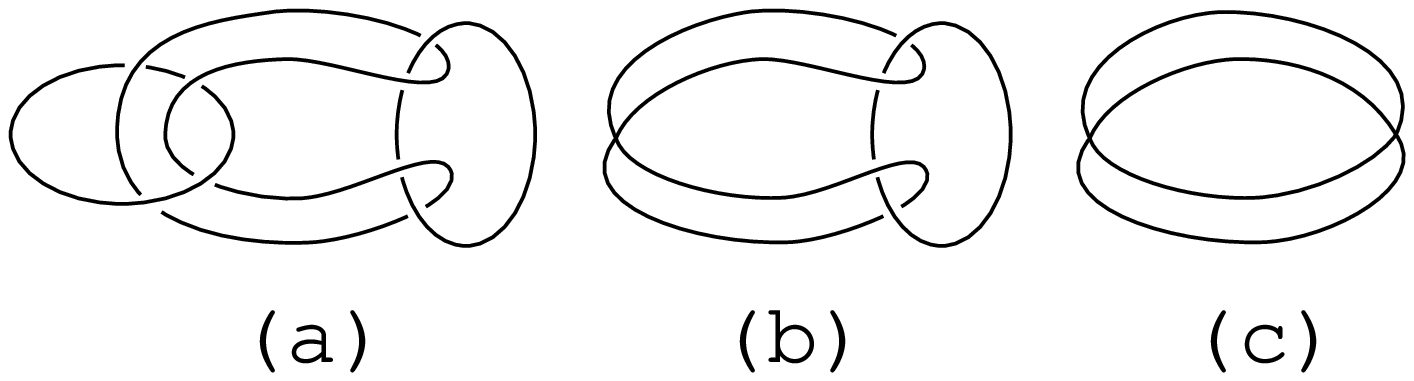}

Fig. 1.2
\end{center}

The graphs $C_n$ with $n>0$ in Fig. 1.3 are shown to be nontrivial in \cite{T-Y} by using Yamada
polynomial
\cite{Yamada}. Note that each $C_n$ is obtained from $C_{n+1}$ by contracting the obvious disk.
Since $C_0$ has no loops and cut-vertices we have that $C_0$ is
irreducible. Therefore we have that each $C_n$ is irreducible. Since a trivial embedding of
a graph abstractly homeomorphic to $C_n$ with $n>0$ has a cutting sphere we have that each
$C_n$ with
$n>0$ is nontrivial.

\begin{center}
\includegraphics[width=5in]{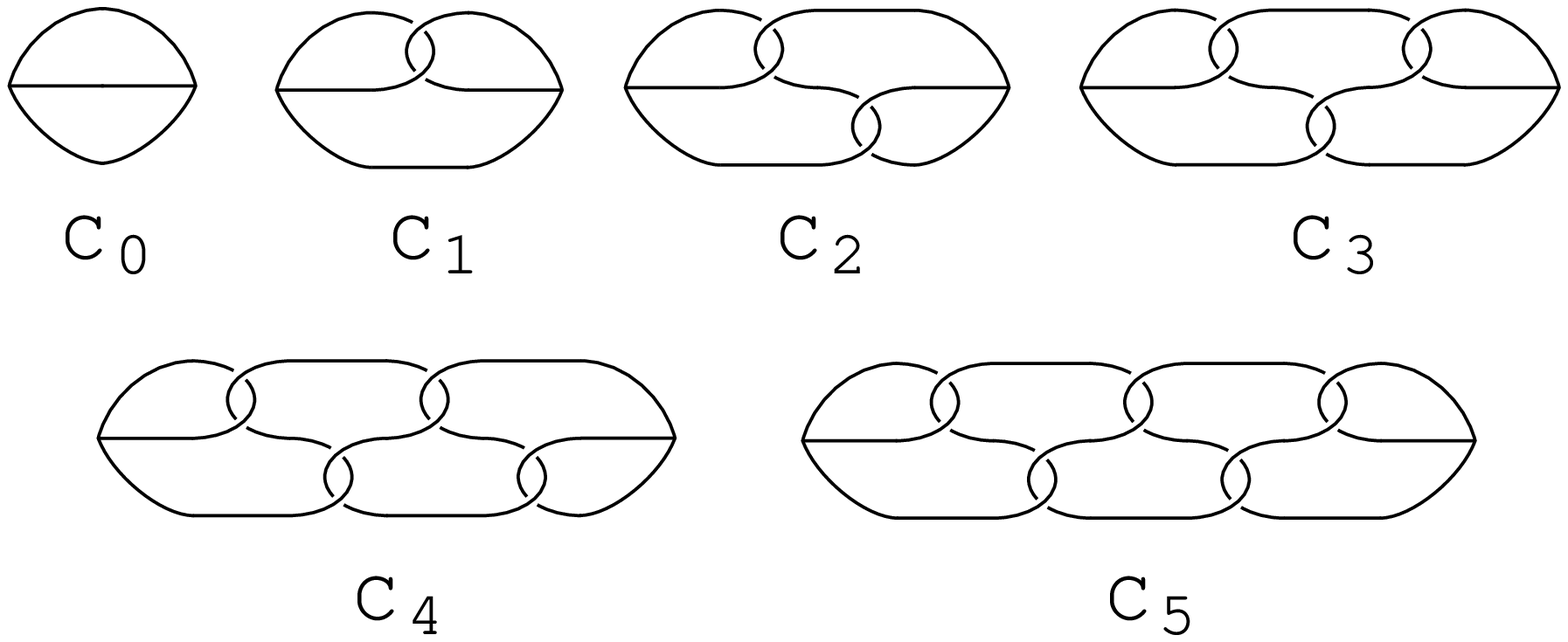}

Fig. 1.3
\end{center}

The graph illustrated in Fig. 1.4 is also irreducible. In fact by contracting the obvious disks we
have a graph without loops and cut-vertices.

\begin{center}
\includegraphics[width=2.5in]{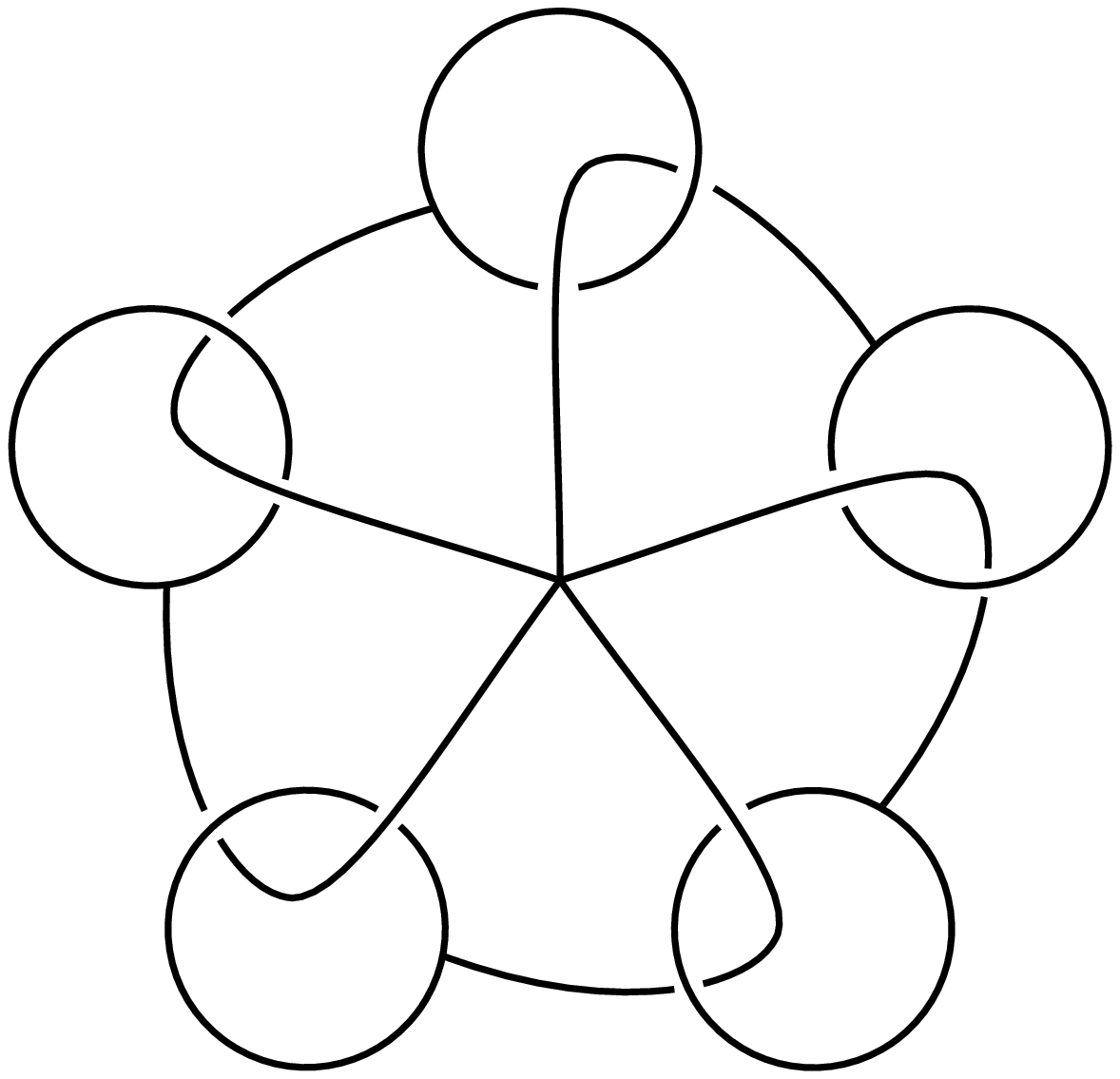}

Fig. 1.4
\end{center}

\vspace{5mm}
\noindent{\large{\bf\S2. Proof}}

\vspace{5mm}\noindent{\bf Proof of Theorem.} It is sufficient to show that if there is a
separating sphere or a cutting sphere of $G$ then there is a separating sphere or a cutting sphere of
$G'$. Let
$S$ be a separating sphere or a cutting sphere of $G$. We may suppose that $S$ and
$D$ are in general position. Then we have that $S\cap D$ consists of some circles and loops. Here
the circles are mutually disjoint and disjoint from the loops. The loops have a point in common.
Namely each loop contains the point $S\cap G$ in case that $S$ is a cutting sphere of $G$.
Therefore in case that $S$ is a separating sphere of $G$ there are no loops.

First suppose that there is a circle or a loop
$\gamma\subset S\cap D$ bounding a disk
$\delta$ in $D$ such that ${\rm int}\delta\cap(S\cup G)=\emptyset$. Then we cut $S$ along $\delta$
and obtain two 2-spheres. It is easy to see that at least one of them is again a separating sphere
or a cutting sphere of $G$. Therefore we may suppose that there are no such circles and loops.

Suppose
that
$S\cap D$ still contains some circles or loops. Let $d$ be a disk in $S$ such that $d\cap
D=\partial d$. Let
$S'$ be a 2-sphere obtained from $d$ by contracting $D$ to a point. It is clear that $S'$
intersects $G'$ at one point. We will show that
$S'$ is good for $G'$. Let $\sigma$ be the disk in $D$ with $\partial\sigma=\partial d$.
Then by the assumption above there is a point
$x\in {\rm int}\sigma\cap G$. Since $D$ is good for $G$ we have that $S'$ is also good for $G'$.
Therefore $S'$ is a cutting sphere of $G'$. Therefore we may suppose that $S\cap D=S\cap
G\neq\emptyset$ or $S\cap D=\emptyset$.

Suppose that $S\cap D=S\cap\partial D=S\cap G\neq\emptyset$. Then by the
assumption on $D$ and the goodness of $S$ for $G$ we have that $S$ becomes a cutting sphere of $G'$
after contracting
$D$ to a point. Suppose that $S\cap D=S\cap{\rm int}D=S\cap G\neq\emptyset$. Then by the goodness of
$S$ for $G$ we have that $S$ becomes a cutting sphere of $G'$ after contracting
$D$ to a point. Suppose that $S\cap D=\emptyset$. Then it is clear that $S$ becomes a separating
sphere or a cutting sphere of $G'$ after contracting
$D$ to a point. This completes the proof. $\Box$

\footnotesize{

\vspace{5mm}

Department of Mathematics, School of Education, Waseda University, 1-6-1 Nishi-Waseda, Shinjuku-ku, Tokyo, 169-8050, Japan

{\it e-mail address}: {\bf taniyama@mn.waseda.ac.jp}

}

\end{document}